\documentclass[a4paper,10pt]{article}
\usepackage[affil-it]{authblk}
\usepackage{amsmath}
\usepackage[utf8x]{inputenc} 
\usepackage{amssymb}
\usepackage{graphicx}
\usepackage{caption} 
\title{ANALYSIS OF A COMBINED NC1-C2 METHOD FOR ELLIPTIC PROBLEM}
\author{Dibyendu Adak%
  \thanks{Electronic address: \texttt{dibyendu.13@iist.ac.in}}}
\affil{Department of Mathematics,
  Indian Institute of Space Science and Technology ,
 Thiruvananthapuram-695547, INDIA}

\author{E. Natarajan%
  \thanks{Electronic address: \texttt{thanndavam@iist.ac.in}; Corresponding author}}
\affil{Department of Mathematics,
  Indian Institute of Space Science and Technology,
 Thiruvananthapuram-695547, INDIA}

\date{Dated: \today}

\begin{document}

\maketitle

\begin{abstract}
It is shown in this paper that non-conforming finite elements on the triangle using $P^{1}$-nonconforming 
polynomials and $P^{2}$ -conforming polynomials can be easily built and used.They  appear as an 'enriched' 
version of the standard piecewise quadratic six-node element.This work is divided into two parts.In the first 
we present the basic- property of the element,namely how it can be built and basic error estimates.We have 
observed that this new element behaves like $P^{1}$ non-conforming element.In the second part we have 
applied our element to the elliptic problem and the theoretical estimate has been guaranteed by numerical result.  
\end{abstract}

\section{INTRODUCTION}
The finite element method has achieved great success in many fields of science and technology since it was 
first suggested in elasticity in the fifth decade of $ 20^{th}$ century.Today it has become a powerful tool for 
solving partial differential equations \cite{brenner2008mathematical,ciarlet2002finite}.The key issue of the finite 
element method is using a discrete solution on the finite element space,usually consisting of piecewise polynomials,to 
approximate the exact solution on the given space according to a certain kind of variational principle.When the 
finite element space is a subspace of the solution space, the method is called conforming.It is known that in this case,the 
finite element solution converges to true solution provided the finite element space approximates the given space in some sense.
   
   In general,  for a $2m$ order elliptic boundary-value problem,the conforming finite element space is a $C^{m-1}$
subspace.It means that the shape function in this conforming finite element space is continuous together with its $m-1$ order derivatives.That is for a second-order problem,  the shape function is continuous and for a fourth order problem,  the shape functions and its derivatives are continuous.It is a rather strong restriction put on the shape functions in the latter case.

   It is prove that to build up a conforming finite element space with $ C^{1}$-continuity for a two dimensional  
fourth-order problem,  like the plate bending problem in elasticity,  at least a quintic polynomial with 18 parameters is required for a triangular element and a bi-cubic polynomial with 16 parameters for a rectangular element.It causes some computational difficulties because the dimension of the related finite element space is fairly large and its structure is rather complicated.

Therefore,  it is desirable to relax directly the $ C^{m-1}$ continuity of the finite element space.It comes to the so called 'nonconforming' finite element method which had and still has a great impact on the development of finite element methods.However,  it was found that some nonconforming element elements converge and some do not. The convergence behaviour sometimes depends on the mesh configuration.

Up to now there has been proposed a vast number of engineer devices based on different mechanical interpolations,like unconventional elements,energy-orthogonal elements with free formulation,quasi-conforming elements,generalised conforming elements and many others.The approximate spaces related to all these elements mentioned above are not included in the given solution space.Hence they are simply called 'nonstandard'     finite elements.A unified  mathematical treatment for analysis of these nonstandard finite elements has been proposed by  various authors.\\

A brief description and analysis of some interesting and important conforming and nonconforming finite elements has been given below.\\

$\textit{$P^{1}$ nonconforming element}$ \cite{ciarlet2002finite,brenner2008mathematical}.  This is a triangular element which is not $C^{0} $. The shape function is a linear polynomial with three nodal parameters at mid points of three edges of the triangle.This element converges for second order elliptic problem with optimal rate.\\

$\textit{Wilson-element} $\cite{taylor1976non,lesaint1976convergence,lesaint1980convergence}.  This is a nonconforming rectangular element.The shape function is a quadratic polynomial with six parameters,  four at vertices of the rectangle and two internal degrees of freedom,  like the second order derivatives.This element converges for rectangular mesh,  but does not converge for arbitrary quadrilaterals.It is interesting to mention that the behaviour of the Wilson element is better than the corresponding bilinear $Q_{1} $ conforming element as many engineering example  have indicated.

$\textit{The rotated $ Q_{1} $ element}$ \cite{li1998nonconforming,ming2001nonconforming,rannacher1990simple}.     $\  $This is a newly established non-conforming rectangular element.The shape function consists of four terms as[1,x,y,$x^{2}-y^{2} $].There are two versions of choosing
nodal parameters.The first one uses four function values at the mid-point of each edge of the rectangle.The second version uses four mean values of the shape function along edges.
Both versions are  convergent for rectangular meshes.However,the first version is not convergent for arbitary quadrilaterals unless certain mess conditions are satisfied.

$\textit{Morley element} $ \cite{morley1968triangular} This is an old and simplest plate element.The shape function is a quadratic polynomial with six nodal parameters.They are three function values at vertices and three normal derivatives at mid-points of three edges.This element does not even belong to $C^{0} $ class,nevertheless it is convergent for the fourth-order problem\cite{lascaux1975some} .  Surprisingly, it is recently proved that the Morley element is divergent for the second-order elliptic problem \cite{nilssen2001robust}.

In contrast, it is well known that there exists for long time 
the $P^{2} $ conforming element  for the second order problem.Its shape function is again a quadratic polynomial with six parameters as three function values at vertices and three function values at mid-points of edges.This quadratic $C^{0}$ element is convergent for the second order problem,  but divergent for the fourth-order problem.

$\textit{Zienkiewicz incomplete cubic triangular element} $  \cite{bazeley1965triangular}.   The shape function consists of incomplete cubic polynomials with specially chosen nine terms.The nine nodal parameters are three function values together with six first partial derivatives with respect to $x$ and $y$ at three vertices.This is a $ C^{0}$ element but not $C^{1}$ .It is proved that this element is convergent only for very special meshes, namely, all edges of triangles are parallel to the three given directions.It is a very interesting phenomenon that the Zienkiewicz
element using the cross-diagonal mesh actually tend to a limit,but it is not the true solution of the given problem,rather of another problem \cite{cai1986limitation} .

Another new nonconforming piecewise quadratic finite element on triangles has been discussed in \cite{fortin1983non}.This element satisfy patch test of Irons and Razzaqque  \cite{irons1972experience}.This implies that on element interfaces,one should ensure the continuity of the approximation at the Gauss-Legendre quadrature points needed for the exact integration of third-degree polynomials. Optimal-error estimates and regularity properties for Dirichlet's problem has been studied by authers.\\
In this paper we will propose a new finite element which is a bridge between conforming and nonconforming finite element.This new element piecewise quadratic and quasi conforming.We have studied error estimation for Dirichlet's problem,and observed that this element does not give optimal convergence rate,which has been generally considered as a major drawback against the use of this element.Though the element is piecewise quadratic but it does not require two point continuity restriction on each interface of $ \tau^{h}$,which is needed for the above mentioned piecewise quadratic nonconforming element \cite{fortin1983non}.We shall show here that these elements are in fact very simple to use and they are nothing but combination of $P^{1} $ nonconforming element with incomplete $P^{2} $ - conforming element.   

\section{Discretization}

\subsection{Continuous Problem}

We consider the following model problem \\
\begin{align}
-\Delta u & =f \quad \text{in} \  \Omega \quad  \\  
        u & =u_{D} \quad \text{on} \  \partial\Omega \quad \nonumber \\  \nonumber
\end{align}
where $\Omega \subset \mathbb{R}^2$ 

\subsection{Notation}
Let $\tau_{h} $ be a conforming triangulation of $ \Omega$.The subscript h refer to the maximum element size $\text h= \text max_{k \in \tau_{h}} \text h_{k} $, where $h_{k}$ is the diameter of an element $K\in \tau_{h} $.$\varepsilon^{h} $ is the set of the edges in $ \tau_{h}$, $\textbf{n} $ is the unit outward normal along $ \partial K$ and the jump $[u] $ across an edge e is a vector defined as follows-\\
Let e be an interior edge shared by two triangles $K_{1}$ and  $K_{2}$ in $ \tau^{h}$ , and $\omega_{j}=\omega \vert_{K_{j}} $ for j=1,2 . We define on e 
\begin{center}
$[\omega]=\omega_{1} \textbf{n}_{1}+\omega_{2} \textbf{n}_{2} $
\end{center}
where $n_{j} $ is the unit normal of e pointing towards the outside of $K_{j} $ .If e is an edge on the boundary of $\Omega $,then we define on e 
\begin{center}
$[\omega]=\omega \textbf{n} $
\end{center}
where $\textbf{n}$ is the unit outer normal of e pointing towards the outside of $\Omega $.


\subsection{Weak Formulation}

Find $u\in V=H^{1}_{0}(\Omega)$ such that 
\begin{center}
$a(u,v)=F(v)\quad \forall \ v \in V=H^{1}_{0}(\Omega) $
\end{center}
 where
\begin{align}
a(u,v) &=\int_{\Omega}\nabla u.\nabla v \\
  F(v) &=\int_{\Omega}fv
\end{align}
$\Omega$ is convex polygon and $f\in L^{2}(\Omega) $ therefore $u\in H^{2}(\Omega) $ by elliptic regularity theory $\cite{s1989topics}.$

\begin{figure}
\caption{NC1-C2 Element}
\centering
\includegraphics[scale=0.6]{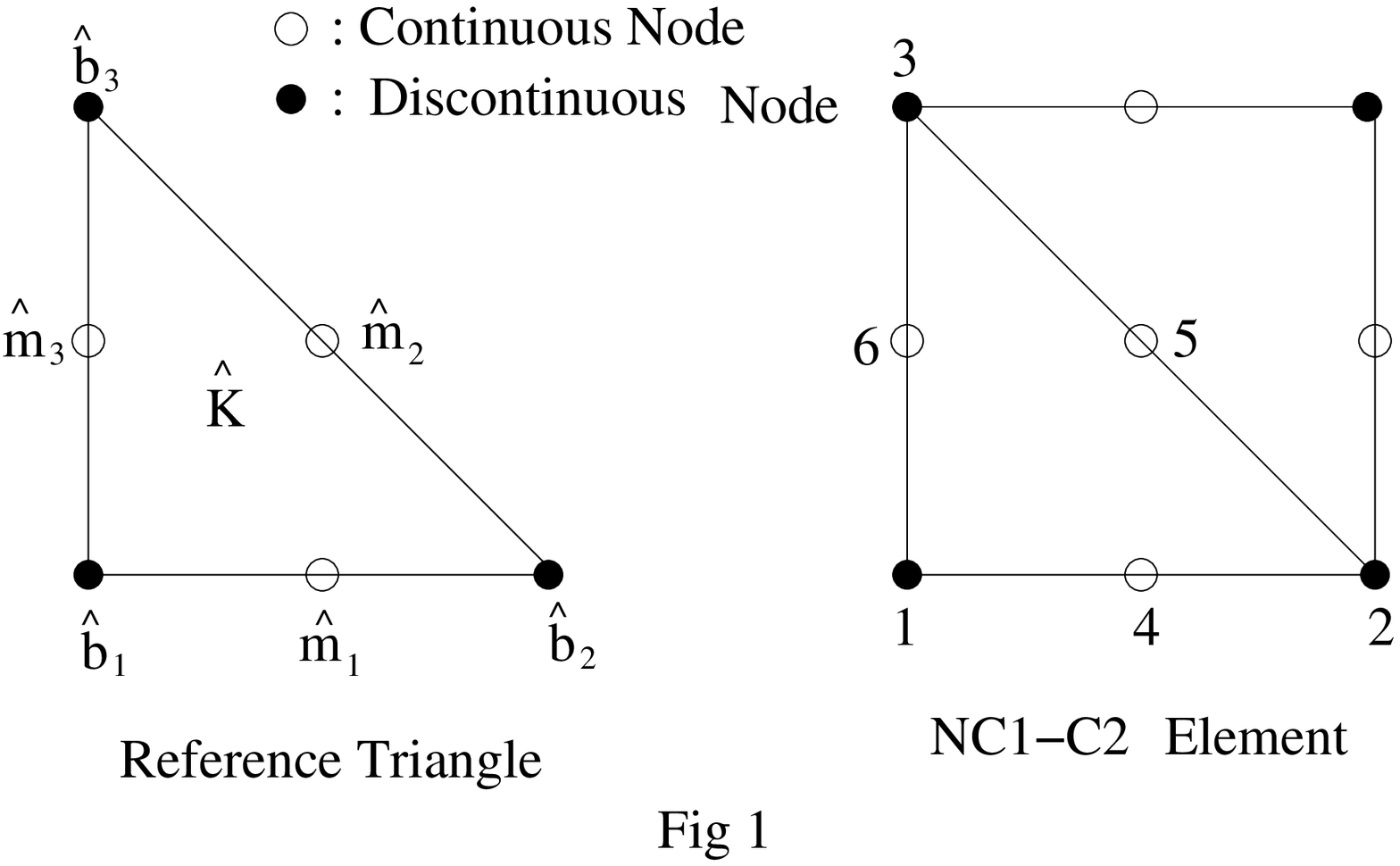}
\end{figure}

$\linebreak $

\subsection{NC1-C2 methods}
In order to define a nonconforming space,we introduce some further notation.\\
$V^{1}_{h}= \{ v \in L^{2}(\Omega) \ : v \vert_{K} \    \text {is linear}\ \forall \  K\in \tau_{h},v \  \text {is continuous at the mid points of the edges of} \  \tau_{h} \} $
The above space is basically $p_{1} $ -nonconforming space.\\
We define
$D_{K}^{2}= \text{span} \{ \phi_{1},\phi_{2},\phi_{3} \} $;where  $\phi_{i}=\hat{\phi}_{i} \circ F^{-1}_{K} $  and $F_{K}$ is Affine mapping from $\hat{K}$ to K [Fig 2].  $\hat{\phi}_{1},\hat{\phi}_{2},\hat{\phi}_{3}$ are basis fn on reference triangle $\hat{K}$ corresponding to the vertices $\hat{b}_{1},\hat{b}_{2},\hat{b}_{3} $ respectively, which is defined by
\begin{align}
\hat{\phi_{1}} &=(-1+2x+2y)(-1+x+y)  \nonumber \\
\hat{\phi_{2}} &=(2x-1)x  \nonumber \\
\hat{\phi_{3}} &=(2y-1)y \nonumber
\end{align}
$\phi_{i}, \  1\leq i \leq 3$ is continuous along edge on each element.In this paper we use the following finite element space-
\begin{center}
 $ V_{h}:= V_{h}^{1} \oplus V_{h}^{2}, \quad  V_{h}^{2}:=  \{v_{h}\in L^{2}(\Omega) : v_{h}\vert_{K} \in D_{K}^{2} \} $ 
 \end{center}

Finite element space $V_{h} $ consists of piecewise quadratic function which is discontinuous along edge of each triangle except at mid points of edges [Fig 1].\\
A typical element $\omega \in V_{h} $ is demonstrated below-
\begin{center}
$\omega=\omega^{1}+\omega^{2} \ \text{where}\  \omega^{1}\in V_{h}^{1}\  \text{and} \ \omega^{2}\in V_{h}^{2} $
\end{center} 
since $\omega^{1}$ is discontinuous along edges except  at the midpoints and $\omega^{2}$ is continuous along edge, $\omega$ is discontinuous along edges except at the mid points.\\
On interior edges jump of typical element is reduced to $P_{1}$ polynomial .Let e be an interior edge which is shared by two triangle $K^{1}$ and $ K^{2}$ . $\omega\vert_{K^{1}} $ and $\omega\vert K^{2} $ are restrictions of $\omega$ on $K^{1}$ and $K^{2}$ respectively.
\begin{equation}
\begin{split}
[\omega] & = \omega\vert K^{1} -\omega \vert K^{2} \\
         & = ( \omega^{1}\vert K^{1}+\omega^{2}\vert K^{1} )-( \omega^{1}\vert K^{2}+ \omega^{2}\vert K^{2})\\
         & = (\omega^{1}\vert K^{1}-\omega^{1}\vert K^{2} )+( \omega^{2}\vert K^{1}-\omega^{2}\vert K^{2})\\
         & = (\omega^{1}\vert K^{1}-\omega^{1}\vert K^{2} )
\end{split}
\end{equation}
This space contains the space of continuous piecewise-quadratic and space of nonconforming piecewise-linear,since
\begin{center}
$V_{h}=V_{h}^{1}+P_{h}^{2}$
\end{center}
where $P_{h}^{2}$ is piecewise-quadratic conforming finite element space.

\subsection{Discretization}
The space $V_{h} $ is not continuous and hence it is no longer in $ H^{1}_{0}(\Omega)$.We must modify the variational form  $a_{h}(.,.) $ in the discretized problem.We define the following bilinear form on $V_{h}+V_{\ast} $.where $V_{\ast} $ is subset of V and exact solution belongs to $V_{\ast} $.   

\begin{align}
a_{h}(v,\omega) &=\sum_{K\in \tau_{h}}\int_{K}\nabla v.\nabla \omega \  dx \\
       F_{h}(\omega)&=\sum_{K\in \tau_{h}} \int_{K}f\omega
\end{align}
\subsection{Consistency}
Let $u \in V_{\ast} \subset V $ such that u satisfies weak formulation,i,e
\begin{equation}
a(u,\omega)=F(\omega)\  \forall\  \omega \in V 
\end{equation}
Then it is obvious that 
\begin{equation}
a_{h}(u,\omega)=F_{h}(\omega) \ \forall \  \omega \in V_{h}
\end{equation}
Hence the discrete bilinear form is consistent.It also satisfies Galerkin's orthogonality condition since
\begin{equation}
a_{h}(u-u_{h},\omega)=0  \ \forall \  \omega \in V_{h}
\end{equation}
\subsection{Discrete stability} The discrete bilinear form $ a_{h}(.,.)  $ enjoys discrete stability on $V_{h} $ if there is $C_{sta}>0 $ such that
\begin{equation}
c_{sta} \parallel v_{h} \parallel \leq \text{sup}_{\omega_{h}\in V_{h}\setminus \{0 \}} \frac{a_{h}(v_{h},\omega_{h})}{\parallel \omega_{h} \parallel}
\end{equation}
\begin{equation}
a_{h}(v_{h},\omega_{h})=\sum_{K\in \tau_{h}}\int_{K}\nabla v_{h}.\nabla \omega_{h}
\end{equation}
considering $v_{h}=\omega_{h}$ \\
\begin{equation}
\begin{split}
\parallel v_{h} \parallel^{2} & =\sum_{K\in \tau_{h}}\int_{K} \vert \nabla v_{h} \vert^{2} \\
                              & =a_{h}(v_{h},v_{h})
\end{split}
\end{equation}
Hence
\begin{equation}
\begin{split}
\parallel v_{h} \parallel & = \frac{a_{h}(v_{h},v_{h})}{\parallel v_{h}\parallel}\\
                          & \leq \text{sup}_{\omega_{h} \in V_{h}\setminus \{ 0\}}\frac{a_{h}(v_{h},\omega_{h})}{\parallel \omega_{h}\parallel}
\end{split}
\end{equation}
Hence $a_{h}(.,.)$ satisfies discrete stability condition.This implies the discrete bilinear form
\begin{equation}
a_{h}(u,\omega)=F_{h}(\omega) \  \forall \ u,\omega \in V_{h}
\end{equation}
is well posed i.e. the discrete bilinear form has unique solution.
\section{Apriori Error Estimation}
$\textbf{Lemma 1}$ $ \ $ Assume dim $V_{h} < \infty $.Let $ a_{h}(..)$ be a symmetric positive definite bilinear form on $ V+V_{h}$ which reduces to a(.,.) on V.Let $u \in V $ solve
\begin{equation}
a(u,v)=F(v) \ \forall \  v\in V \nonumber
\end{equation}
where $F\in V' \cap V_{h}' $.Let $u_{h} \in V_{h} $ solve 
\begin{equation}
a_{h}(u_{h},v)=F(v) \  \forall \  v \in V \nonumber
\end{equation}
Then 
\begin{equation}
\parallel u-u_{h} \parallel_{h} \leq \text{inf}_{v \in V_{h}} \parallel u-v \parallel _{h}+\text{sup}_{\omega \in V_{h}\setminus \{ 0\}} \frac{\vert a_{h}(u-u_{h},\omega)}{\parallel \omega \parallel_{h}}
\end{equation}
where $\parallel . \parallel_{h}=\sqrt{a_{h}(.,.)} $\\

proof:- Let $\tilde{u_{h}} \in V_{h} $ satisfies 
\begin{equation}
a_{h}(\tilde{u_{h}},v)=a_{h}(u,v) \  \forall \  v \in V_{h} \nonumber
\end{equation}
which implies that 
\begin{equation}
a_{h}(\tilde{u_{h}}-u,v) =0 \ \forall \ v\in V_{h} \nonumber
\end{equation}
 \begin{equation} 
\Rightarrow \ \parallel u-\tilde{u_{h}} \parallel_{h} =\text{inf}_{v \in V_{h}} \parallel u-v \parallel_{h} \nonumber
\end{equation}
Then
\begin{equation}
\begin{split}
\parallel u-u_{h}\parallel_{h} &\leq \parallel u- \tilde{u_{h}} \parallel_{h}+\parallel \tilde{u_{h}}-u_{h} \parallel_{h} \\                
                               &\leq \parallel u-\tilde{u_{h}} \parallel_{h}+\text{sup}_{\omega \in V_{h}\setminus \{0 \}} \frac{\vert a_{h}(\tilde{u_{h}}-u_{h},\omega)\vert}{\parallel \omega \parallel_{h}} \nonumber
\end{split}
\end{equation}
$\linebreak $
$\textbf{Lemma 2} $ $\quad$ Let K be an arbitrary element of conforming triangulation $\tau_{h}$.Then the following inequality holds
\begin{equation}
\vert e \vert^{-1} \parallel \zeta \parallel_{L^{2}(e)}^{2} \leq \  C(h_{K}^{-2} \parallel \zeta  \parallel_{L^{2}(K)}^{2}+ \vert \zeta \vert_{H^{1}(K)}^{2})\ \forall \  \zeta \in H^{1}(K)
\end{equation}
where $\vert e \vert $ denotes the length of edge $e\subset \partial K$, $h_{K}$ =diam$\ $K,and the positive constant depends only on the chunkiness parameter of K.\\

Proof: See the details in \cite{brenner2008mathematical}\\

$\textbf{Lemma 3} $ \quad  Let all  assumptions of Lemma 2 hold and $\omega $ be an arbitrary element of $ V_{h} $.Then

\begin{equation}
\vert e \vert  \parallel [\omega] \parallel_{L^{2}(e)}^{2} \leq \ C \  \sum_{K \in \tau_{e}} h_{K}^{2} \vert \omega \vert_{H^{1}(K)}^{2} 
\end{equation}
where $[\omega] $ denotes jump of $\omega $ along  edge  $e \in \varepsilon^{h}$.\\

Proof: Using lemma 2 we can write
\begin{align}
 \vert e \vert^{-1} \parallel [\omega] \parallel_{L^{2}(e)}^{2} \ &\leq \  C \sum_{K \in \tau_{e}}(h_{K}^{-2} \parallel \omega \nonumber \parallel_{L^{2}(K)}^{2} +\vert \omega \vert_{H^{1}(K)}^{2} ) \\ \nonumber
 \vert e \vert \parallel [\omega] \parallel_{L^{2}(e)}^{2} &\leq \ C \sum_{K \in \tau_{e}} (\vert e \vert^{2} h_{K}^{-2} \parallel \omega \parallel_{L^{2}(K)}^{2}+\vert e \vert^{2} \vert \omega \vert_{H^{1}(K)}^{2}) \nonumber
\end{align}
Where $\tau_{e} $ is the set of triangles in $ \tau_{h}$ containing e on their boundaries.
Again $[\omega]$ =0 at midpoint of each edge of K hence first part of of right-hand side will be vanished.Therefore we have 
\begin{equation}
\begin{split}
\vert e \vert \parallel [\omega] \parallel_{L^{2}(e)}^{2} &\leq \ C \sum_{K \in \tau_{e}} \vert e \vert^{2} \vert \omega \vert_{H^{1}(K)}^{2}\\ \nonumber
                                                               & \leq C \sum_{K \in \tau_{e}}  h_{K}^{2} \vert \omega \vert_{H^{1}(K)}^{2} \nonumber
\end{split}
\end{equation}
$ \linebreak $
$\textbf{Lemma 4}$ \quad Let B be a ball in $\Omega $ such that $\Omega  $ is star-shaped with respect to B and such that its radius $ \rho > (\frac{1}{2}) \rho_{max}$. Let $Q^{m}u $ be the Taylor polynomial of order m of u averaged over B where $u \in W_{p}^{m}(\Omega) $ and $p\geq 1$.Then 
\begin{equation}
\vert u-Q^{m}u \vert_{W^{k}_{p}(\Omega)} \leq C_{m,n,\gamma} d^{m-k} \vert u \vert_{W_{p}^{m}(\Omega)} \  k=0,1,\dot{...} ,m, 
\end{equation}
where d=diam$(\Omega)$ and $ \rho_{max}$=sup $\{ $ $\rho:\Omega$ is star-shaped with respect to a ball of radius $\rho$  $\}$

Proof: See the details in \cite{brenner2008mathematical}\\
\subsection{}

The important ingredient in the error analysis is a bound on the approximation error $\parallel u-u_{I} \parallel $ where $u_{I} \in V_{h} $ is a suitable interpolation  which agrees with u at mid point of each edges of $\varepsilon^{h}$ of exact solution u.The interpolation operator is defined at the element level. We just require the local approximation property
\begin{equation}
\vert u - u_{I} \vert_{H^{s}(K)} \leq C h^{p+1-s}_{K} \vert u \vert_{H^{p+1}(K)} \quad \forall \  K \in \tau_{h} \ ,s=0,1,2 \nonumber
\end{equation} 
It will be useful to define it explicitly.It is defined in two steps.\\
\begin{figure}
\caption{Affine Mapping}
\centering
\includegraphics[scale=0.5]{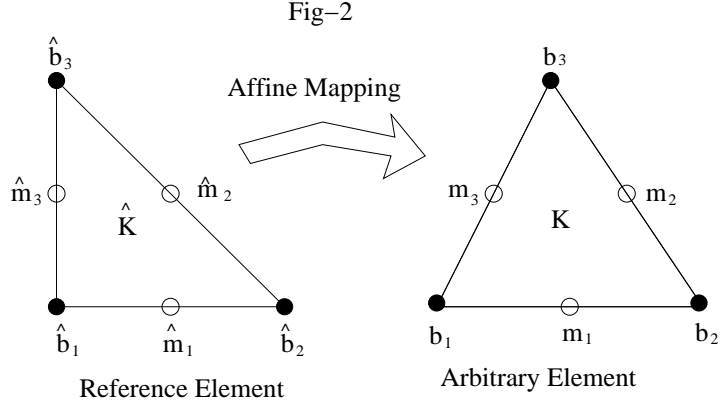}
\end{figure}
Let $\hat{K} $ be the reference triangle with vertices$ \  \hat{b_{1}},\hat{b_{2}},\hat{b_{3}} \  $  whose coordinates are $\  $ (0,0),(1,0),(0,1)  respectively and $\hat{m_{i}} $ be the midpoint of the side joining i and i+1(modulo 3) vertices.\\
We define 
\begin{equation}
\hat{I}^{1}(\hat{u})=\hat{u}(\hat{m}_{1}) \hat{\phi}_{4}+\hat{u}(\hat{m}_{2})\hat{\phi}_{5}+\hat{u}(\hat{m}_{3})\hat{\phi}_{6} \nonumber
\end{equation}
\begin{equation}
\hat{I}^{2}(\hat{u})=(\hat{u}(\hat{b}_{1})-\hat{I}^{1}(\hat{u})(\hat{b}_{1}))\hat{\phi}_{1}+(\hat{u}(\hat{b}_{2})-\hat{I}^{1}(\hat{u})(\hat{b}_{2}))\hat{\phi}_{2}+(\hat{u}(\hat{b}_{3})-\hat{I}^{1}(\hat{u})(\hat{b}_{3}))\hat{\phi}_{3}  \nonumber
\end{equation}
Finally we define interpolation as
\begin{equation}
\begin{split}
\hat{I}(\hat{u}) & = \hat{I}^{1}(\hat{u})+ \hat{I}^{2}(\hat{u}) \\
                 & = \sum_{j=1}^{6} \hat{L}_{j}(\hat{u})\hat{\phi}_{j}
\end{split}
\end{equation}
where $\hat{L}_{i} $ for i=1,2,3,4,5,6 $\ $ continuous linear functional.\\
$\linebreak$
Now we will show that $P_{2}(\hat{K}) $ is unisolvent with respect to these functionals, i.e.for an arbitrary polynomial $\hat{p} \in P_{2}(\hat{K}) $ $ \ $ $ \hat{L}_{i}(\hat{p})=0$ implies $\hat{p}=0 $. \\
Since $\hat{p} \in P_{2}(\hat{K}) $ implies that $\hat{p} $ can be written as linear combination of basis of $P_{2}(\hat{K}) $ .Then 
\begin{align}
  \  \hat{p} \ &=\sum_{i=1}^{6} C_{i} \phi_{i}. \nonumber \\  \nonumber
 C_{4} & =\hat{p}(\hat{m}_{1})=\hat{L}_{4}(\hat{p})=0\\ \nonumber
C_{5} & =\hat{p}(\hat{m}_{2})=\hat{L}_{5}(\hat{p})=0\\  \nonumber
 C_{6} & =\hat{p}(\hat{m}_{3})=\hat{L}_{6}(\hat{p})=0\\  \nonumber
 \Rightarrow \hat{p} &=C_{1}\hat{\phi}_{1}+C_{2}\hat{\phi}_{2}+C_{3}\hat{\phi}_{3}
 \end{align}
 again
 \begin{align}
  \ C_{1} &=\hat{p}(\hat{b}_{1})=\hat{L}_{1}(\hat{p})=0  \nonumber \\ \nonumber
 C_{2} & =\hat{p}(\hat{m}_{2})=\hat{L}_{2}(\hat{p})=0\\ \nonumber
 C_{3} & =\hat{p}(\hat{m}_{3})=\hat{L}_{3}(\hat{p})=0\\  \nonumber
 \Rightarrow \ \hat{p} &=0
\end{align}
Similarly it can be shown that for arbitrary $ \hat{p}\in P_{2}(\hat{K})$,$\ $ $\hat{I}(\hat{p})=\hat{p} $ .         
 Let $(K,P_{2}(K),\Sigma)$ be an affine finite element of $(\hat{K},P_{2}(\hat{K}),\hat{\Sigma)}$ $\ $ where$\ $ $\Sigma=\{L_{1},L_{2},L_{3},L_{4},L_{5},L_{6}\} $ and $ \ $  $\hat{\Sigma}=\{\hat{L}_{1},\hat{L}_{2},\hat{L}_{3},\hat{L}_{4},\hat{L}_{5},\hat{L}_{6}\} $ .\\
Then for all $v \in H^{m+1}(K) $ we have 
\begin{equation}
\parallel D^{s}(u-I_{h}u)\parallel_{L^{2}(K)} \leq C\ h_{K}^{m+1-s}\ \parallel D^{m+1} u \parallel_{L^{2}(K)}
\end{equation}
where $s \leq m+1 $ and m=0,1,2.
 
\subsection{}
We want to find out $\parallel u-u_{h}\parallel_{h}$
Using lemma 1 we can write
\begin{equation}
\parallel u-u_{h}\parallel_{h} \leq \text{inf}_{v \in V_{h}} \parallel u-v \parallel_{h}+\text{sup}_{\omega \in V_{h}\setminus \{ 0\}} \frac{a_{h}(u-u_{h},\omega)}{\parallel \omega \parallel_{h}} \nonumber
\end{equation}
\begin{equation}
\begin{split}
\text{inf}_{v \in V_{h}} \parallel u-v \parallel_{h} &\leq \parallel u-u_{I} \parallel_{h}\\
                                                     &\leq \  C\  h^{2}\  \vert u \vert_{H^{3}(\Omega)} \quad \text{by} \quad (20) 
\end{split}
\end{equation}
Now we have to estimate 
\begin{center}
$\text{sup}_{\omega \in V_{h}\setminus \{ 0\}} \frac{\vert a_{h}(u-u_{h},\omega)\vert}{\parallel \omega \parallel_{h}}$
\end{center}
\begin{equation}
\begin{split}
a_{h}(u-u_{h},\omega) & = \sum_{K \in \tau_{h}} \int_{K} \nabla u. \nabla \omega dx -\int_{\Omega}f \omega dx\\
                      &=\sum_{k\in \tau_{h}} [\int_{\partial k} \nabla u.\omega \textbf{n} ds - \int_{K} \Delta u \omega d\Omega]-\int_{K}f\omega d \Omega \\
                      &=\sum_{K \in \tau_{h}} \int_{\partial K } \nabla u.\omega \textbf{n} ds\\
                      &=\sum_{e \in \varepsilon^{h}} \int_{e} \nabla u.[\omega] ds 
\end{split}
\end{equation}
Again we have the following estimate
\begin{equation}
\sum_{e\in \varepsilon^{h}} \int_{e} \nabla u.[\omega] ds=\sum_{e \in \varepsilon^{h}} \int_{e}(\nabla u.\textbf{n}_{e}-c_{e})\textbf{n}_{e}.[\omega]ds
\end{equation}
since $ \int_{e}[\omega]=0$ \quad by (4).
\begin{equation}
\begin{split}
\sum_{e \in \varepsilon^{h}} \int_{e} \nabla u.[\omega] ds &=\sum_{e \in \varepsilon ^{h}} \int_{e} (\nabla u.\textbf{n}_{e}-c_{e})\textbf{n}_{e}.[\omega]ds \\
                                                           & \leq \sum_{e \in \varepsilon^{h}} \vert e \vert^{-\frac{1}{2}} \parallel \nabla u. \textbf{n}_{e} -c_{e} \parallel_{L^{2}(e)} \vert e \vert^{\frac{1}{2}} \parallel [\omega] \parallel_{L^{2}(e)}\\
                                                           & \leq (\sum_{e \in \varepsilon^{h}} \vert e \vert^{-1} \parallel \nabla u.\textbf{n}_{e}-C_{e} \parallel^{2}_{L^{2}(e)} )^{\frac{1}{2}}(\sum_{e \in \varepsilon^{h}} \vert e \vert \parallel[\omega] \parallel^{2}_{L^{2}(e)})^{\frac{1}{2}}\\
                                                           & \leq C [\sum_{e \in \varepsilon^{h}} \text{min}_{K \in \tau_{e}}(h_{K}^{-2} \parallel \nabla u .\textbf{n}_{e}-C_{e} \parallel^{2}_{L^{2}(T)}+\vert u \vert ^{2}_{H^{2}(K)})]^{\frac{1}{2}} [\sum_{e \in \varepsilon^{h}} \sum_{K \in \tau_{e}} h_{K}^{2} \vert \omega \vert^{2}_{H^{1}(K)}  ] ^{\frac{1}{2}} \  \text{by} \ (16,17) \\                                                         
                                                           & \leq C h \vert u \vert _{H^{2}(\Omega)} \parallel \omega \parallel_{h} \nonumber
\end{split}
\end{equation}
Hence we have
\begin{equation}
\vert a_{h}(u-u_{h},\omega)\vert  \leq C\  h \  \vert u \vert_{H^{2}(\Omega)} \parallel \omega \parallel_{h} \nonumber 
\end{equation}
\begin{equation}
\Rightarrow \quad \frac{\vert a_{h}(u-u_{h},\omega)\vert}{\parallel \omega \parallel_{h}}   \leq C \ h \ \vert u \vert_{H^{2}(\Omega)} \nonumber
\end{equation}
\begin{equation}
\Rightarrow \quad \parallel u-u_{h} \parallel_{h}  \leq C\ h \ \vert u \vert_{H^{2}(\Omega)}\quad \text{using \ (21)} 
 \end{equation}

\subsection{$ L^{2} $-Error Estimate}
Let $ \eta \in H^{2}(\Omega)\bigcap H^{1}_{0}(\Omega) $ satisfy
\begin{equation}
a(\eta,v)=\int_{\Omega} v (u-u_{h}) \ \forall \ v \in H^{1}_{0}(\Omega) \nonumber
\end{equation}
and $\eta_{h} \in V_{h} $ satisfy
\begin{equation}
a_{h}(\eta_{h},v)=\int_{\Omega} v \ (u-u_{h})  \  \forall \ v \in V_{h}  \nonumber
\end{equation}
Again we have 
\begin{equation}
\begin{split}
\parallel u-u_{h} \parallel^{2}_{L^{2}(\Omega)} &= \int_{\Omega}(u-u_{h})^{2} \\
                                                &= \int_{\Omega}(u-u_{h})(u-u_{h}) \\
                                                &=\int_{\Omega}u(u-u_{h})-\int_{\Omega}u_{h}(u-u_{h})\\
                                                &=a(u,\eta)-a_{h}(u_{h},\eta_{h}) \\
                                                &=a_{h}(u,\eta)-a_{h}(u_{h},\eta)+a_{h}(u_{h},\eta)-a_{h}(u_{h},\eta_{h})\\
                                                &=a_{h}(u-u_{h},\eta)+a_{h}(u_{h},\eta-\eta_{h})\\
                                                &=a_{h}(u-u_{h},\eta-\eta_{h})+a_{h}(u-u_{h},\eta_{h})+a_{h}(u_{h},\eta-\eta_{h}) 
\end{split}
\end{equation}

Using  estimation (24) we can write
\begin{equation}
\begin{split}
a_{h}(u-u_{h},\eta-\eta_{h}) & \leq \parallel u-u_{h} \parallel_{h} \parallel \eta-\eta_{h} \parallel_{h} \quad (\text{using \ Cauchy-Schwarz  }) \\
                              & \leq \  C \  h^{2} \  \vert u \vert_{H^{2}(\Omega)} \ \vert \eta \vert_{H^{2}(\Omega)} 
\end{split}
\end{equation}
Where $\ $ C $\  $ is generic constant.
\begin{equation}
a_{h}(u-u_{h},\eta_{h})=a_{h}(u-u_{h},\eta_{h}-\eta_{I})+a_{h}(u-u_{h},\eta_{h}) \nonumber
\end{equation} 

\begin{equation}
\begin{split}
a_{h}(u-u_{h},\eta_{h}-\eta_{I}) & \leq \parallel u-u_{h} \parallel_{h} \  \parallel \eta_{h} -\eta_{I} \parallel_{h} \\
                                 & \leq \  C \ h^{2} \ \vert u \vert_{H^{2}(\Omega)} \ \vert \eta  \vert_{H^{2}(\Omega)} \nonumber
\end{split}
\end{equation}

\begin{equation}
\begin{split}
\parallel \eta_{h}-\eta_{I}\parallel_{h} & = \parallel \eta_{h}-\eta+\eta-\eta_{I} \parallel_{h} \\
                                         & \leq \parallel \eta_{h}-\eta \parallel_{h} +\parallel \eta - \eta_{I} \parallel_{h} \\
                                         &\leq \ C_{1}\ h \ \vert \eta \vert_{H^{2}(\Omega)} + C_{2} \  h \ \vert \eta \vert_{H^{2}(\Omega)} \\
                                         &=\ C \ h \ \vert \eta \vert_{H^{2}(\Omega)} 	\nonumber
\end{split}
\end{equation}
\begin{equation}
\begin{split}
a_{h}(u-u_{h},\eta_{I}) &= \sum_{e \in \varepsilon^{h}} \  \int_{e}\ \nabla u.[\eta_{I}] \ ds \quad \quad by \ (22)  \\
                        & =\sum_{e \in \varepsilon^{h}} \  \int_{e}\ \nabla u.[\eta_{I}-\eta] \ ds \\
                        & =\sum_{e \in \varepsilon^{h}} \ \int_{e} \ (\nabla u.\textbf{n}_{e}-C_{e})\ \textbf{n}_{e}.[\eta_{I}-\eta]\ ds  \nonumber
\end{split}
\end{equation}
Where $ \   C_{e} \ $ is arbitrary constant and we have used the following two identity
\begin{center}
$ [\eta]=0 \quad \& \quad \int_{e}[\eta_{I}]=0   $
\end{center}

\begin{equation}
\begin{split} 
a_{h}(u-u_{h},\eta_{I}) & \leq \sum_{e \in \varepsilon^{h}} \vert e \vert^{-\frac{1}{2}} \parallel \nabla u.\textbf{n}_{e}-C_{e} \parallel_{L^{2}(e)} \vert e \vert^{\frac{1}{2}} \parallel \eta-\eta_{I} \parallel_{L^{2}(e)} \nonumber \\ 
                        & \leq C \ [\sum_{e \in \varepsilon^{h}} \text{min}_{K \in \tau_{e}}(h^{-2}_{K} \text{inf}_{C_{e} \in \mathbb{R}} \parallel \nabla u.\textbf{n}_{e} -C_{e} \parallel^{2}_{L^{2}(K)} + \vert u \vert^{2}_{H^{2}(K)} )]^{\frac{1}{2}} \\
                        & \times [\sum_{e \in \varepsilon^{h}} \vert e \vert^{2} \sum_{K \in \tau_{e}} (h^{-2}_{K} \parallel \eta-\eta_{I} \parallel^{2}_{L^{2}(K)} + \vert \eta-\eta_{I}\vert^{2}_{H^{1}(K)} )]^\frac{1}{2} \\
                        &\leq\  C \ \vert u \vert_{H^{2}(\Omega)} (h^{2} \ \vert \eta \vert_{H^{2}(\Omega)})\quad \text{using} \ (18,20) 
\end{split}
\end{equation}
Therefore
\begin{equation}
a_{h}(u-u_{h},\eta_{h}) \leq \ C \ h^{2} \ \vert u \vert_{H^{2}(\Omega)} \  \vert \eta \vert_{H^{2}(\Omega)}
\end{equation}
Again we have
\begin{equation}
\begin{split}
a_{h}(u_{h},\eta-\eta_{h}) & = a_{h}(\eta-\eta_{h},u_{h}) \\
                           & =a_{h}(\eta-\eta_{h},u_{h}-u_{I})+a_{h}(\eta-\eta_{h},u_{I}) 
\end{split} 
\end{equation}

\begin{equation}
a_{h}(\eta-\eta_{h},u_{h}-u_{I}) \  \leq \ \parallel \eta-\eta_{h} \parallel_{h} \parallel u_{h}-u_{I} \parallel_{h}  \nonumber
\end{equation}
and 
\begin{equation}
\begin{split}
\parallel u_{h}-u_{I}\parallel_{h} &= \parallel u_{h}-u+u-u_{I}  \parallel_{h} \\
                                   & \leq \ \parallel u-u_{h} \parallel_{h}+\parallel u-u_{I} \parallel_{h} \\
                                   & \leq \ C \ h \vert u \vert_{H^{2}(\Omega)}  \nonumber
\end{split}
\end{equation}

Therefore
\begin{equation}
a_{h}(\eta-\eta_{h},u_{h}-u_{I}) \leq \ C \ h^{2} \vert \eta \vert_{H^{2}(\Omega)} \  \vert u \vert_{H^{2}(\Omega)}
\end{equation}
We can write second term of (28) as 
\begin{equation}
\begin{split}
a_{h}(\eta-\eta_{h},u_{I}) &= \sum_{e \in \varepsilon^{h}} \int_{e} \nabla \eta.[u_{I}] \\
                           &= \sum_{e \in \varepsilon^{h}} \int_{e} (\nabla \eta.\textbf{n}_{e}-C_{e})\ \textbf{n}_{e}.[u-u_{I}] \ ds  \nonumber
\end{split}
\end{equation}
where we have used the following two identity
\begin{center}
$[u]=0 \quad \& \quad \int_{e}[u_{I}]= 0 $ ,e is an interior edge.
\end{center}
Therefore
\begin{equation}
\begin{split}
a_{h}(\eta-\eta_{h},u_{I}) &= \sum_{e \in \varepsilon^{h}} \ \int_{e} (\nabla \eta.\textbf{n}_{e}-C_{e})\ \textbf{n}_{e}.[u-u_{I}] \ ds \\
                           & \leq \sum_{e \in \varepsilon^{h}} \ \vert e \vert^{-\frac{1}{2}} \parallel \nabla \eta .\textbf{n}_{e}-C_{e} \parallel_{L^{2}(e)} \ \vert e \vert^{\frac{1}{2}} \parallel u-u_{I}\parallel_{L^{2}(e)} \\
                           & \leq \ C \ \vert \eta \vert_{H^{2}(\Omega)}(h^{2} \ \vert u \vert_{H^{2}(\Omega)} ) 
\end{split}
\end{equation}
Hence
\begin{equation}
\begin{split}
a_{h}(u_{h},\eta-\eta_{h}) &=a_{h}(\eta-\eta_{h},u_{h}-u_{I})+a_{h}(\eta-\eta_{h},u_{I}) \\
                           & \leq \ C \ h^{2} \ \vert \eta \vert_{H^{2}(\Omega)} \ \vert u \vert_{H^{2}(\Omega)}  
\end{split}
\end{equation}
Using  estimations (26),(27),(31) in (25) we can write
\begin{equation}
\parallel u-u_{h}\parallel^{2}_{L^{2}(\Omega)}  \leq \ C \ h^{2} \ \vert \eta \vert_{H^{2}(\Omega)}\ \vert u \vert_{H^{2}(\Omega)} 
\end{equation}
using elliptic regularity we can write
\begin{equation}
\parallel \eta \parallel_{H^{2}(\Omega)}  \leq \ C \ \parallel u-u_{h} \parallel_{L^{2}(\Omega)} 
\end{equation}
Therefore
\begin{equation}
 \quad \parallel u-u_{h} \parallel_{L^{2}(\Omega)} \leq \ C \ h^{2} \  \vert u \vert_{H^{2}(\Omega)} \nonumber
\end{equation} 
 
\section{Numerical Test}
In this section we perform grid convergence studies for the proposed 'NC1-C2' method.The new method was implemented  using penalization technique.We estimate the experimental order of convergence by the formula
\begin{equation}
Eoc=log(\frac{E(2h)}{E(h)}) \nonumber
\end{equation} 
where $E(h)=\parallel u-u_{h}\parallel $ is the error in the specified norm.The result indicate the same convergence behaviour for NC1-C2 method and and $p^{1}$-nonconforming method.We have glued second order conforming  element with first order nonconforming element suitably and observed that first order nonconforming element dominated second order conforming  element and we have obtained an EOC of 2.0 in the $L^{2}$ norm and $1.0 $ in the $H^{1} $ norm which is same as $p^{1}$ nonconforming method.EOC stands for experimental order of convergence.\\ 

\textbf{Stationary Diffusion Problem :}
 We consider the Poisson equation 
\begin{equation}
-\Delta u=f \quad \text{on} \quad \Omega=(0,1)\times(0,1)
\end{equation} 
with homogeneous boundary conditions and the right hand side $\ f=2\pi^{2}sin(\pi x)sin(\pi y)$. The exact solution is given by $u(x,y)=sin(\pi x)sin(\pi y) $.\\

\section*{Acknowledgments}
Authors thank Prof. Alexei Lozinski(Université de Franche-Comté, Besancon) and Prof Frederic Hecht(LJLL-Université Pierre et Marie Curie) 
for fruitful discussions on implementation of NC1-C2 method during the workshop CIMPA-2015 held at IIT BOMBAY, INDIA.
 
\bibliographystyle{unsrt}
\bibliography{arch1}

\begin{thebibliography}{10}

\bibitem{brenner2008mathematical}
Susanne~C Brenner and Ridgway Scott.
\newblock {\em The mathematical theory of finite element methods}, volume~15.
\newblock Springer Science \& Business Media, 2008.

\bibitem{ciarlet2002finite}
Philippe~G Ciarlet.
\newblock {\em The finite element method for elliptic problems}, volume~40.
\newblock Siam, 2002.

\bibitem{taylor1976non}
Robert~L Taylor, Peter~J Beresford, and Edward~L Wilson.
\newblock A non-conforming element for stress analysis.
\newblock {\em International Journal for Numerical Methods in Engineering},
  10(6):1211--1219, 1976.

\bibitem{lesaint1976convergence}
Pierre Lesaint.
\newblock On the convergence of wilson's nonconforming element for solving the
  elastic problems.
\newblock {\em Computer Methods in Applied Mechanics and Engineering},
  7(1):1--16, 1976.

\bibitem{lesaint1980convergence}
P~Lesaint and M~Zlamal.
\newblock Convergence of the nonconforming wilson element for arbitrary
  quadrilateral meshes.
\newblock {\em Numerische Mathematik}, 36(1):33--52, 1980.

\bibitem{li1998nonconforming}
Bo~Li and Mitchell Luskin.
\newblock Nonconforming finite element approximation of crystalline
  microstructure.
\newblock {\em Mathematics of Computation of the American Mathematical
  Society}, 67(223):917--946, 1998.

\bibitem{ming2001nonconforming}
Pingbing Ming and Zhong-Ci Shi.
\newblock Nonconforming rotated element for reissner--mindlin plate.
\newblock {\em Mathematical Models and Methods in Applied Sciences},
  11(08):1311--1342, 2001.

\bibitem{rannacher1990simple}
Rolf Rannacher, Stefan Turek, and Deutsche Forschungsgemeinschaft.
\newblock {\em A simple nonconforming quadrilateral Stokes element}.
\newblock Universit{\"a}t Heidelberg. SFB 123, 1990.

\bibitem{morley1968triangular}
LSD Morley.
\newblock The triangular equilibrium element in the solution of plate bending
  problems.
\newblock {\em Aero. Quart}, 19:149--169, 1968.

\bibitem{lascaux1975some}
P~Lascaux and P~Lesaint.
\newblock Some nonconforming finite elements for the plate bending problem.
\newblock {\em Revue fran{\c{c}}aise d'automatique, informatique, recherche
  op{\'e}rationnelle. Analyse num{\'e}rique}, 9(1):9--53, 1975.

\bibitem{nilssen2001robust}
Trygve Nilssen, Xue-Cheng Tai, and Ragnar Winther.
\newblock A robust nonconforming 𝐻$^2$-element.
\newblock {\em Mathematics of Computation}, 70(234):489--505, 2001.

\bibitem{bazeley1965triangular}
GP~Bazeley, Yo~K Cheung, Bo~M Irons, and OC~Zienkiewicz.
\newblock Triangular elements in plate bending—conforming and nonconforming
  solutions.
\newblock In {\em Proceedings of the Conference on Matrix Methods in Structural
  Mechanics}, pages 547--576. Wright Patterson AF Base, Ohio, 1965.

\bibitem{cai1986limitation}
W~Cai.
\newblock The limitation problem for zienkiewicz triangle elements.
\newblock {\em Math. Numer. Sinica}, 8:345--353, 1986.

\bibitem{fortin1983non}
M~Fortin and M~Soulie.
\newblock A non-conforming piecewise quadratic finite element on triangles.
\newblock {\em International Journal for Numerical Methods in Engineering},
  19(4):505--520, 1983.

\bibitem{irons1972experience}
Bruce~M Irons and Abdur Razzaque.
\newblock Experience with the patch test for convergence of finite elements.
\newblock {\em The mathematical foundations of the finite element method with
  applications to partial differential equations}, 557:587, 1972.

\bibitem{s1989topics}
S.~Kesavan.
\newblock {\em Topics in functional analysis and applications}.
\newblock John Wiley \& Sons, 1989.

\end{thebibliography}

\newpage
\begin{figure}[H]
\includegraphics[width=12cm]{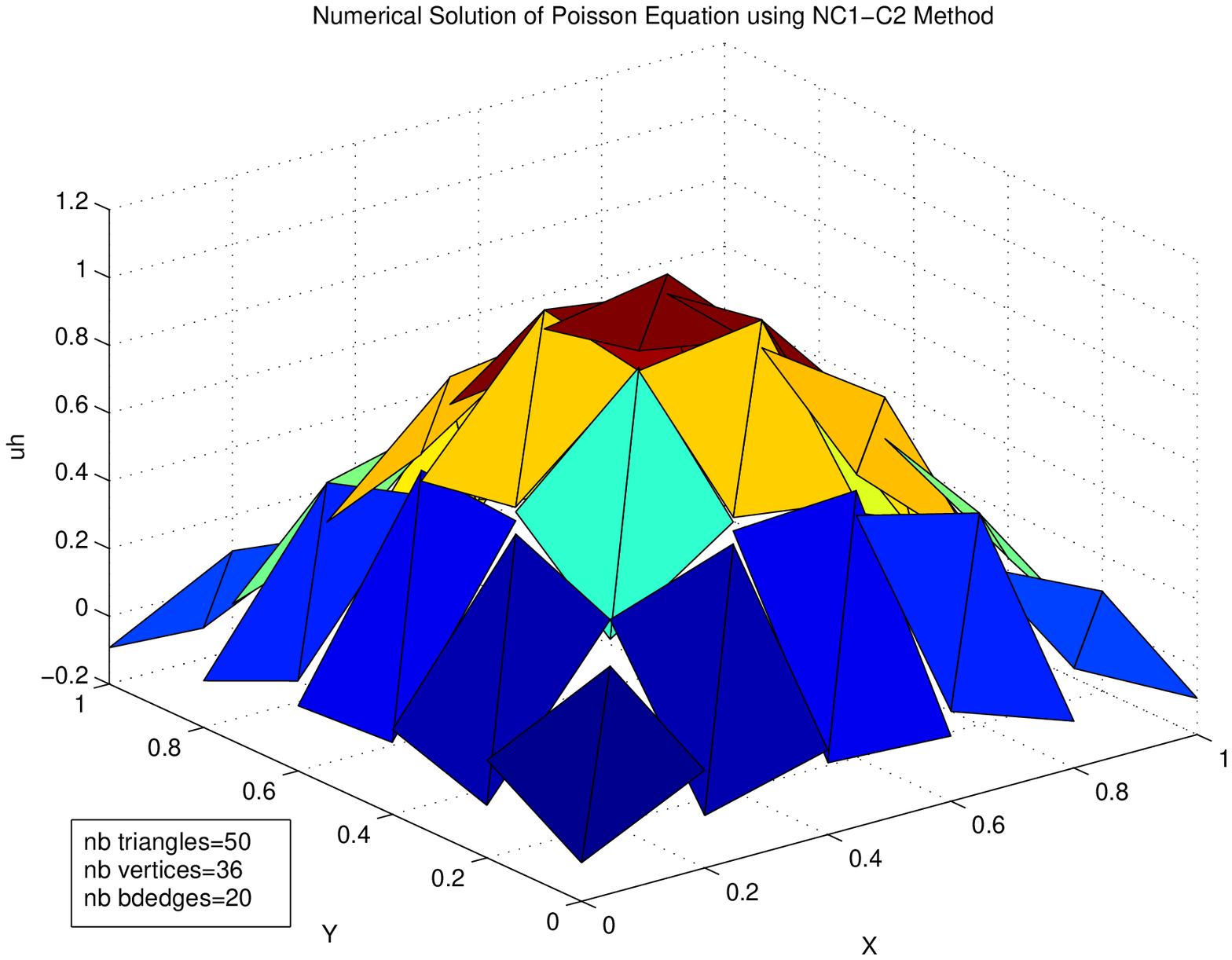}
\end{figure}

\newpage
\begin{figure}[H]
\includegraphics[width=12cm]{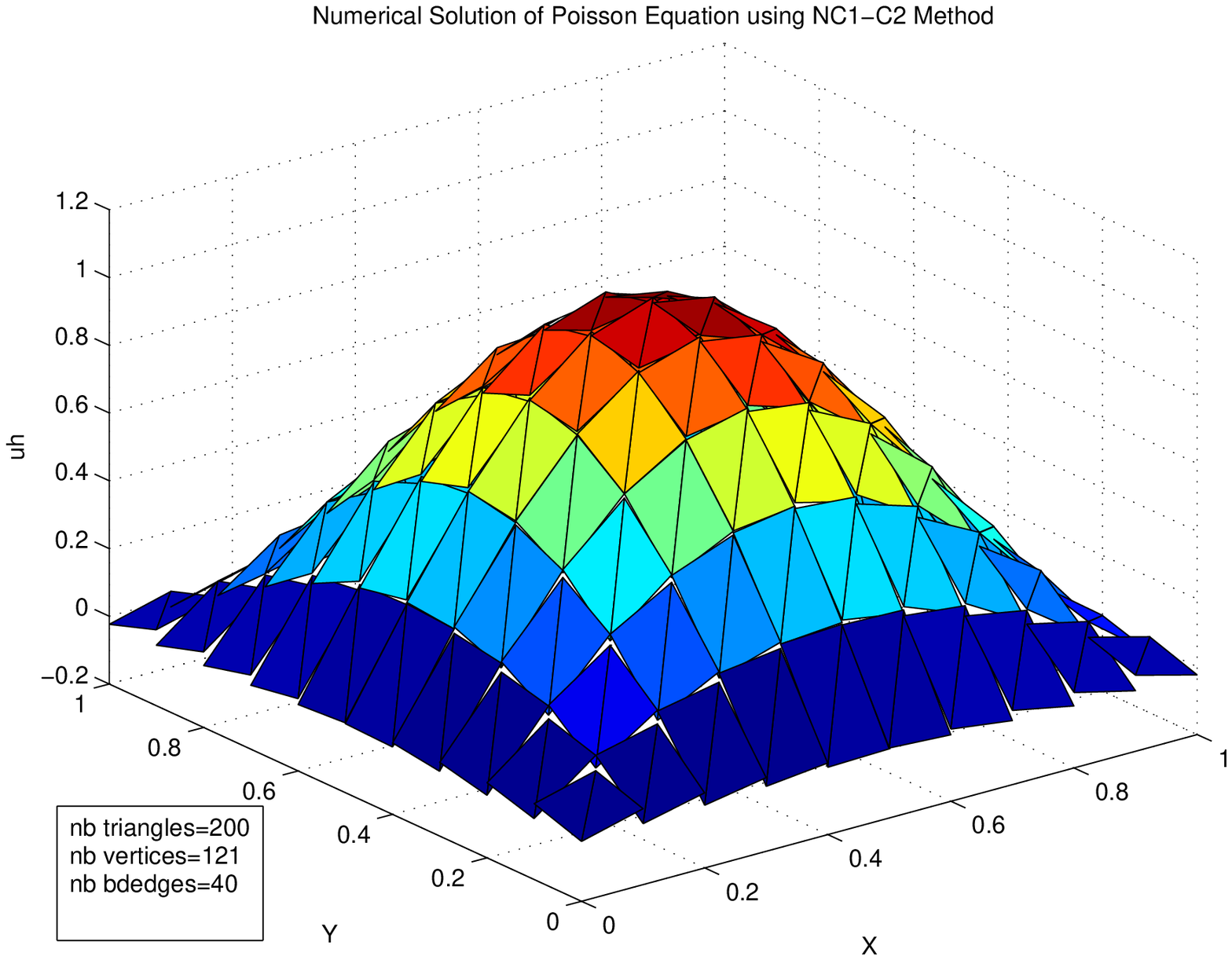}
\end{figure}

\begin{figure}[H]
\includegraphics[width=12cm]{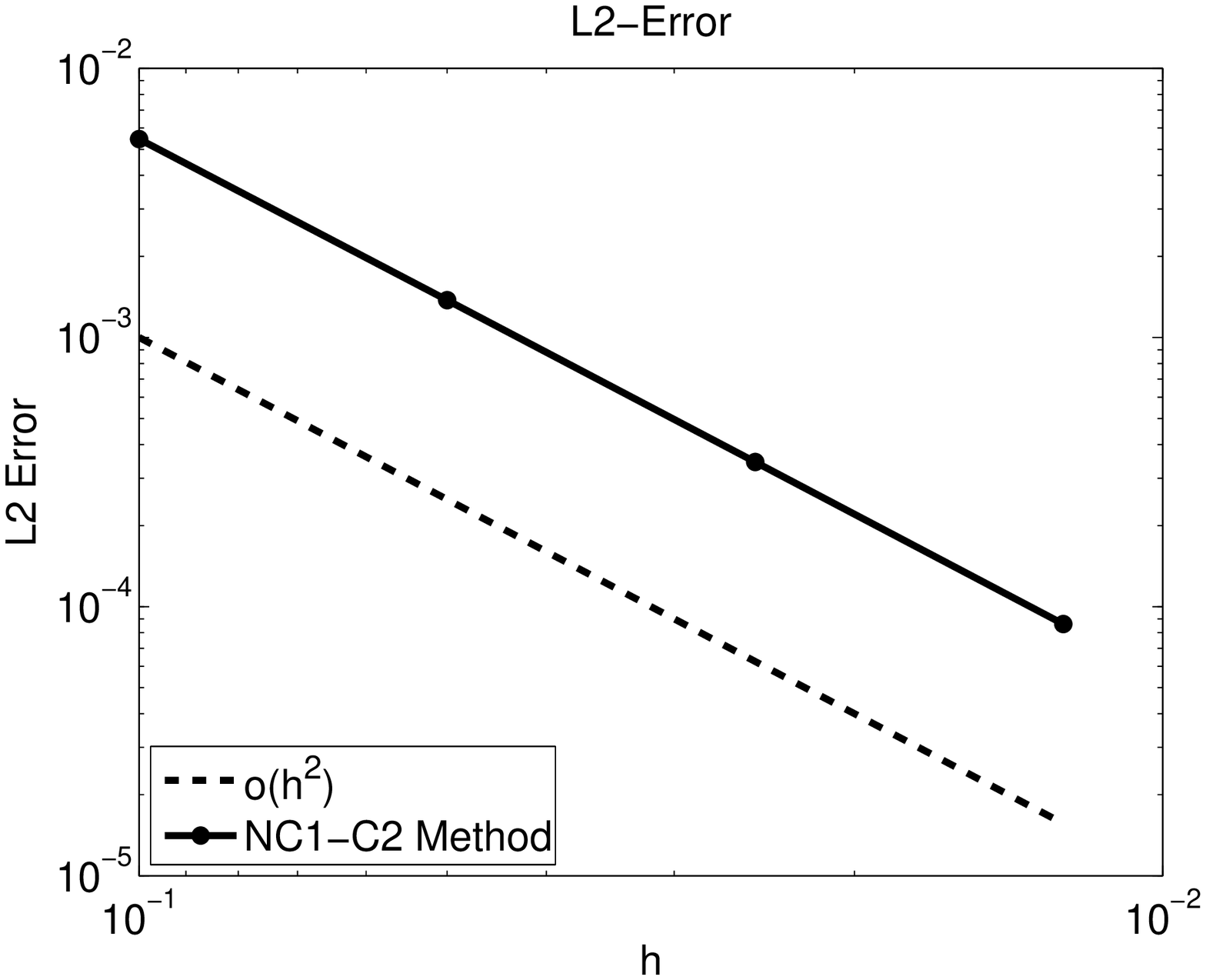}
\end{figure}

\newpage
\begin{figure}[H]
\includegraphics[width=12cm]{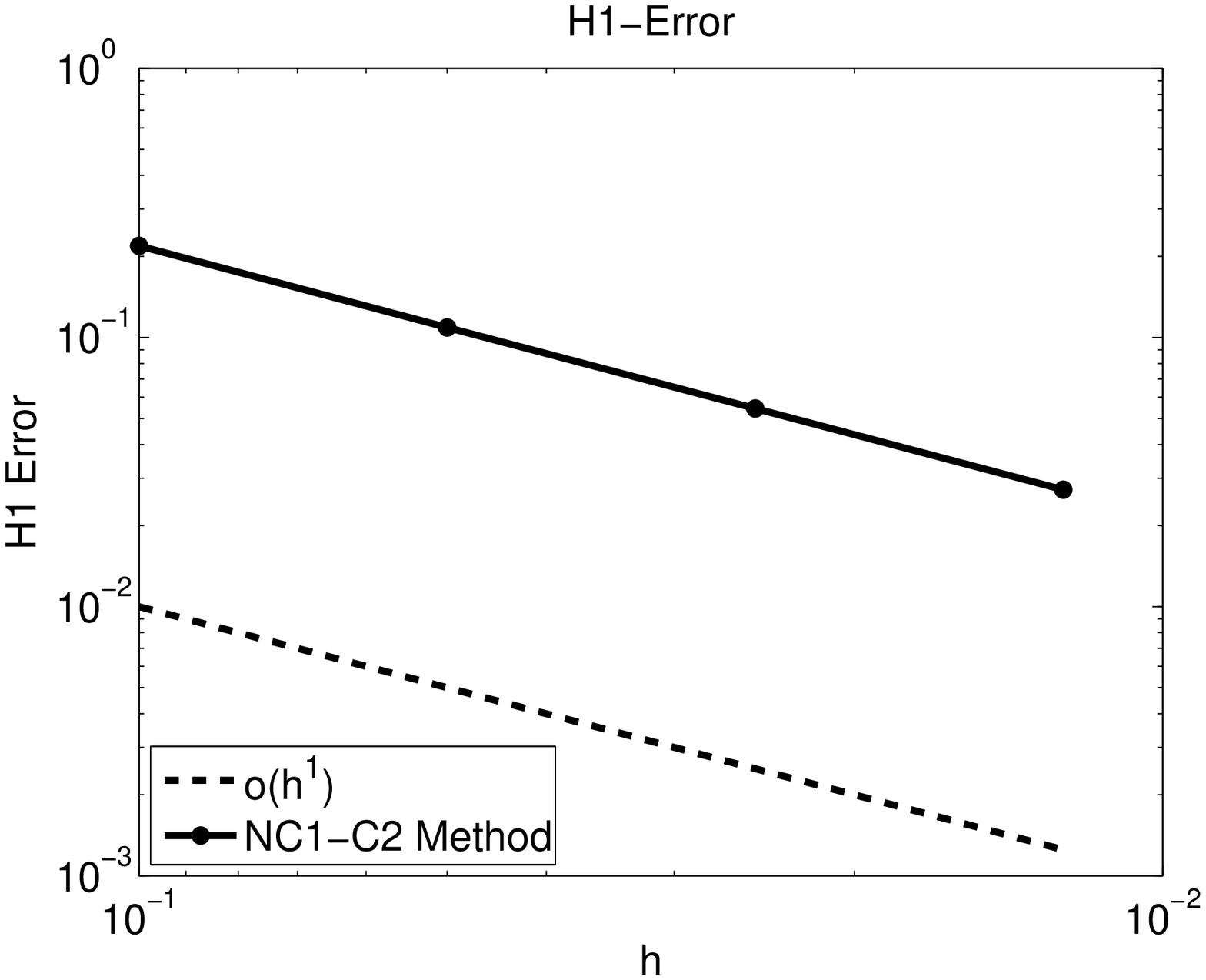}
\end{figure}
\end{document}